\documentclass[12pt]{amsart}

\usepackage{amssymb}
\usepackage{graphicx}
\usepackage{epsf}
\usepackage{amsmath}
\usepackage[latin1]{inputenc}
\usepackage{amsfonts}
\usepackage{color}

\newtheorem{prop}{Proposition}

\DeclareMathOperator{\Arg}{Arg}
\DeclareMathOperator{\II}{Im}

\def\sgn{{\rm sgn}}

\def\bz{{\bf z}}

\def\t{_\theta}

\def\C{{\mathbb C}}

\def\Z{{\mathbb Z}}

\def\Real{{\mathbb R}}

\def\U{{\mathcal U}}

\def\In{{\rm In}}
\def\Out{{\rm Out}}

\title[A symmetry property of some harmonic algebraic curves]{A symmetry property\\ of some harmonic algebraic curves}

\author[Aval, Marckert]{Jean-Christophe Aval, Jean-Fran\c{c}ois Marckert}\address[J.-C. Aval, J.-F. Marckert]{Universit\'e de Bordeaux \\ LaBRI, CNRS\\ 351 cours
 de la Lib\'eration\\ 33405 Talence cedex\\ FRANCE}
\email{aval@labri.fr, marckert@labri.fr}
\urladdr{http://www.labri.fr/perso/aval, http://www.labri.fr/perso/marckert}

\thanks{This work has been supported by the ANR project MARS (BLAN06-2$\_$0193)} 

\begin{document} 

\maketitle
\begin{abstract}
The aim of this note is to give a surprising symmetry property of some harmonic algebraic curves: when all the roots $z_i$ of a complex polynomial $P$ lie on the unit circle $\U$, the points of $\U$ different from the $z_i$, and such that $\Arg(P(z))=\theta$, form a regular $n$-gon, where $n$ is the degree of $P$.
\end{abstract}      

\medskip

Let  $\bz=\{z_1,\dots,z_n\}$ be a multiset of $n$ points in the complex plane $\C$ and $P$ the monic polynomial with root set $\bz$:
\[P(z)=\prod_{i=1}^n (z-z_i).\]
For $\theta$ a fixed real number of your choice, consider
$$C\t(P)=\{z\in \C\ :\ \II(e^{-i\theta} P(z))=0\}.$$
The set $C\t(P)$ coincides up to $\bz$, to the set $\{z \in \C :\ \Arg(P(z))=\theta [\pi]\}$.
These curves arise in the Gauss approach to the Fundamental Theorem of Algebra (see e.g. Stillwell \cite{Still}, and Martin \& al. \cite{MSS}). In their paper Martin \& al. \cite{MSS} and then Savitt \cite{savitt} initiated the study of the combinatorial topology of the families $C\t(P)$. The idea are the following ones: the curves $C\t(P)$ have $2n$ asymptotes at angles $(\pi k +\theta)/n$, for $k\in\{0,\dots,2n-1\}$, and form in the generic case $n$ non intersecting curves.  This induces a matching: $k$ and $k'$ are matched if and only if the asymptotes $(\pi k+\theta)/n$ and $(\pi k'+\theta)/n$ lie on the same connected component in $C\t(P)$. The papers \cite{MSS} and \cite{savitt} aim at studying these matchings, and also the properties of the so-called necklaces, formed by the families of matchings obtained when $\theta$ traverses the set $[0,\pi]$. \par
\begin{figure}[htbp]
\centerline{\includegraphics[height= 3cm]{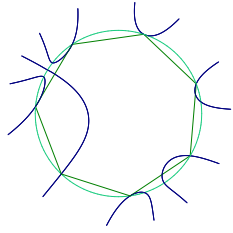}}
\caption{\label{figu} An example where $n=7$, $\theta=0$ and the roots $z_i$ randomly chosen.}
\end{figure}

Let us now state and prove our result. The set $\bz$ is clearly included in $C\t(P)$. It turns out that when $\bz$ is included in the unit circle $\U=\left\{z~:~|z|=1\right\}$, the set $C\t(P)\cap \U$ presents a quite surprising symmetry -- illustrated at Figure \ref{figu} -- that can be stated as follows.
\begin{prop}\label{poor-form}If $\bz$ is a subset of $\U$, then 
\[C\t(P)\cap \U=\bz\cup G(\bz)\]where $G(\bz)$  is the regular $n$-gon on $\U$, with set of vertices $\Big\{e^{i(\Omega +2 k \pi/n)}, k=1,\dots,n\Big\},$ for
$$\Omega:=\frac{2\theta-\sum_{j=1}^n \Arg(z_j)}{n}-\pi.$$
\end{prop}
There exists a purely geometric proof of this Proposition using that the measure of a central angle is twice that of the inscribed angle intercepting the same arc; we provide below a more compact analytic proof. 
\proof 
We will only consider $z\notin \bz$. We have the equivalence:
$$z\in C\t(P)\setminus \bz\ \ \Longleftrightarrow\ \ z\notin \bz,~\sum_{i=1}^n \Arg(z-z_i)=\theta\ \ [\pi],$$
where $\Arg(z)\in\Real/2\pi\Z$ stands for (any chosen determination of) the argument of $z\neq0$. 
Now for any $\nu$ and $\psi$ real numbers,
$$e^{i\nu}-e^{i\psi}=e^{i\frac{\nu+\psi}2}\big(e^{i\frac{\nu-\psi}2}-e^{i\frac{-\nu+\psi}2}\big)=2i\sin ((\nu-\psi)/2)e^{i\frac{\nu+\psi}2}.$$
Thus
$$\Arg(e^{i\nu}-e^{i\psi})=\frac{\nu+\psi}2 +\frac{\pi}2+\pi\times\sgn\big(\sin ((\nu-\psi)/2)\big)\ \ [2\pi]$$ 
Hence, $z\in C\t(P)\setminus\bz$ is equivalent to: 
$$ z\notin \bz,~~\sum_{j=1}^n \Big(\frac{\Arg(z)+\Arg(z_j)}2+\frac{\pi}2\Big)=\theta\ \ [\pi], $$
which leads to the conclusion at once.~$\Box$
%Since we are working modulo $\pi$, adding or substracting $\pi/2$ is identical, the sign of the sinus does not matter. It follows that $z\in C\t(P)\setminus\bz$ if and only if:
%$$z\notin \bz,~~\Arg(z)=\frac{2\theta-\sum_{j=1}^n \Arg(z_j)}{n}-\pi\ \ \Big[\frac{2\pi}n\Big].~~~~\Box$$

\noindent\bf Note. \rm If $z_i$ is a root of multiplicity $k$ of $P$, and if $z_i$ belongs to $G(\bz)$, then in the neighborhood of $z_i$, $C\t(P)$ has $k$ tangents, one of them coinciding with the tangent of the circle at $z_i$. Moreover, it is simple to check that if $z_i$ is not on $G(\bz)$, then the tangents of $C\t(P)$ at $z_i$ are not tangent to $\U$.

%%%%%%%%%%%%%%%%%%%%%%%%%%%%%%

\end{document}